\documentclass[12pt]{amsart}
\usepackage{amsmath}
\usepackage{amsfonts}
\usepackage{amssymb}
\usepackage{bbding}
\usepackage{txfonts}
\usepackage[shortlabels]{enumitem}
\usepackage{ifpdf}
\ifpdf
  \usepackage[final,hyperindex]{hyperref}
\else
  \usepackage[final,hyperindex,hypertex]{hyperref}
\fi

\newtheorem{thm}{Theorem}[section]%
\newtheorem{lem}[thm]{Lemma}%
\newtheorem{prop}[thm]{Proposition}%

\theoremstyle{definition}%
\newtheorem{defn}[thm]{Definition}%
\newtheorem{exa}[thm]{Example}%
\begin{document}
\title[Monomial Rota-Baxter operators]{Classification of monomial Rota-Baxter operators on $\mathbf{k}[x]$}

\author{Houyi Yu}
\address{
School  of Mathematics and Statistics, Southwest University, Chongqing, 400715, China}
\vspace{-30mm}
\address{Mathematics Department, Polytechnic School of Engineering, New York University \\
New York, $11201$, USA}
\email{yuhouyi@swu.edu.cn}

%========================================================================
\hyphenpenalty=8000
\date{\today}

\begin{abstract}
Rota-Baxter operators were introduced to solve certain analytic and combinatorial
problems and then applied to many fields in mathematics and mathematical
physics. The polynomial algebra  $\mathbf{k}[x]$ plays a central role both in analysis and algebra.
In this paper, we explicitly classified  all monomial Rota-Baxter operators
on $\mathbf{k}[x]$.
\end{abstract}



\maketitle

%\tableofcontents

\hyphenpenalty=8000 \setcounter{section}{0}

%========================================================================

\section{Introduction}\label{Introduction}
Rota-Baxter operator is an algebraic abstraction and generalization of the integration by parts formula of calculus.
The study originated from the probability study of Baxter \cite{Baxter1960} in 1960
and then developed by the school of Rota \cite{Rota1995}.
This concept has  been closely related to many fields in mathematics and mathematical physics
such as combinatorics,
number theory, operads , quantum field theory
(see \cite{{ConnesKreimer1998},guokeigher2000,{Ebrahimiguokreimer2002},Ebrahimiguokreimer2004,Bai2007,guozhang2008,Baiguoni2010}
and the references therein). See \cite{guo2009} for a brief introduction and \cite{guo2012} for a more detailed treatment.

Because of the importance of Rota-Baxter operators,
it is useful to explicitly determine their classification.
In recent years, some progress
regarding computations of Rota-Baxter operators on semigroup algebras and Lie algebras have been achieved,
with applications to pre-Lie algebras, dendriform algebras and the classical Yang-Baxter equation
\cite{lxb2008,peibaifuo2014,GLBJ2014,rbopaiaoguo2014}.

The polynomial algebra $\mathbf{k}[x]$ is an important object both in analysis and in algebra. It provides an ideal
testing ground to see how an abstractly defined Rota-Baxter operator is related to the integration
operator, because of its analytic connection, as functions, and its
algebraic significance as a free object in the category of $\mathbf{k}$-algebras.
Guo, Rosenkranz and Zheng \cite{rbopaiaoguo} studied Rota-Baxter operators on the polynomial algebra  $\mathbf{k}[x]$ that send monomials to monomials
and give a sufficient condition for a monomial linear
operator on $k[x]$ to be a Rota-Baxter operator.

In this paper we further investigate the detailed calssification of monomial Rota-Baxter operators on $\mathbf{k}[x]$.
In Section \ref{sectionweightzero} we first give a necessary and sufficient condition for a monomial linear operator to be a Rota-Baxter operator
of weight zero by improving \cite[Theorem 3.3]{rbopaiaoguo}, and then give a specific construction for this kind of operators.
Section \ref{sectionweightnonzero} is devoted to the monomial Rota-Baxter operators  of weight nonzero.
We show that this kind of operators can be divided into four classes.

\section{Monomial Rota-Baxter operators of weight zero}\label{sectionweightzero}

We first recall some background and preliminary results that will be used in this paper.

Throughout the paper, unless otherwise stated, we assume that $\mathbf{k}$ is an integral domain containing the field $\mathbb{Q}$ of rational numbers,
the monoid of natural numbers (nonnegative integers) is denoted by $\mathbb{N}$,
we write $\mathbb{N}^*$ for the semigroup of positive integers.

\begin{defn}
Let $\mathbf{k}$ be a commutative ring and let $\lambda$ be an element of $\mathbf{k}$.
A {\bf Rota-Baxter operator of weight $\lambda$} on a commutative $\mathbf{k}$-algebra $R$ is defined to be a $\mathbf{k}$-linear
map $P: R\rightarrow R$ that satisfies the {\bf Rota-Baxter equation}
\begin{equation}\label{rtequ}
P(x)P(y)=P(xP(y))+P(P(x)y)+\lambda P(xy) \quad for\ all\quad x,y\in R.
\end{equation}
\end{defn}

\begin{defn}\label{defnpbttht}
A linear operator $P$ on $\mathbf{k}[x]$ is called {\bf monomial} if for each $n \in \mathbb{N}$, we have
\begin{equation*}%\label{defnofp}
P(x^n)=\beta(n)x^{\theta(n)}\qquad {\rm with}\qquad \beta: \mathbb{N}\rightarrow \mathbf{k}
\qquad{\rm and}\qquad  \theta: \mathbb{N}\rightarrow \mathbb{N}.
\end{equation*}
If $\beta(n)=0$, then the value of $\theta(n)$ does not matter; by convention we set $\theta(n)=0$ in this case.
A monomial operator is called {\bf degenerate} if $\beta(n) = 0$ for some $n\in \mathbb{N}$.
\end{defn}

Let $A$ be a nonempty set and let $B$ be a set containing a distinguished element $0$. For a
map $\phi: A \rightarrow B$ we define  $\mathcal{Z}_\phi:=\{a\in A|\phi(a)=0\}$ to be the  {\bf zero set} of $\phi$.
Then we write its  {\bf support} as $\mathcal{S}_\phi:= A\backslash\mathcal{Z}_\phi$.
Thus a monomial linear operator $P$ on $\mathbf{k}[x]$ is nondegenerate
if and only if $\mathcal{Z}_\beta= \emptyset$.
By Definition \ref{defnpbttht}, we have $\mathcal{Z}_\beta\subseteq\mathcal{Z}_\theta$, so that $\mathcal{S}_\theta\subseteq\mathcal{S}_\beta$.
A straightforward calculation (see \cite[Lemma 3.2]{rbopaiaoguo}) shows that
 $\mathcal{S}_\beta=\mathcal{S}_\theta$ and
$\mathcal{Z}_\beta=\mathcal{Z}_\theta$ for a monomial Rota-Baxter operator $P$ of weight zero.
However, it is possible even if $\mathcal{S}_\beta\cap\mathcal{S}_\theta=\emptyset$ for a monomial Rota-Baxter operator of
weight nonzero as shown in Example \ref{examp21}.

In this section, all Rota-Baxter operators are assumed to be of weight $\lambda=0$ defined by $P(x^n)=\beta(n)x^{\theta(n)}$, $n\in \mathbb{N}$.
We will give a specific classification for this kind of Rota-Baxter operators.
We first give a necessary and sufficient condition for monomial Rota-Baxter operators by improving \cite[Theorem 3.3]{rbopaiaoguo} as follows.
\begin{thm}\label{equicodmonomial}
Let $P$ be a monomial linear operator on $\mathbf{k}[x]$ defined by $P(x^n)=\beta(n)x^{\theta(n)}$, $n\in \mathbb{N}$.
Then $P$ is a Rota-Baxter operator of weight zero if and only if $\theta$ and $\beta$ satisfy the following conditions:
\begin{enumerate}
\item\label{bicontain1} $\mathcal{Z}_\beta+\theta(\mathcal{S}_\beta)\subseteq \mathcal{Z}_\beta$,
\ $\mathcal{S}_\beta+\theta(\mathcal{S}_\beta)\subseteq \mathcal{S}_\beta$;

\item\label{bicontain2}  for all $m,n\in \mathcal{S}_\beta$, we have
\begin{align}
\theta(m)+\theta(n)=\theta(\theta(m)+&n)=\theta(m+\theta(n)),\label{thetaeq3}\\
\beta(m)\beta(n)=\beta(m)\beta(\theta(m)+&n)+\beta(n)\beta(m+\theta(n)).\label{betaeq4}
\end{align}
\end{enumerate}
\end{thm}
\begin{proof}
In view of \cite[Theorem 3.3]{rbopaiaoguo}, we only need to show the fact that  $P$ is a Rota-Baxter operator  of weight zero implies   $\mathcal{S}_\beta+\theta(\mathcal{S}_\beta)\subseteq \mathcal{S}_\beta$, which follows from Lemma \ref{sbetasubset1}.
\end{proof}

\begin{lem}\label{m+thetam}
Let $P$ be a monomial Rota-Baxter operator on $\mathbf{k}[x]$.
Then for any $m\in \mathcal{S}_\beta$ and any nonnegative integer $k$, $m+k\theta(m)\in \mathcal{S}_\beta$.
Furthermore,
\begin{align}
\theta(m+k\theta(m))&=(k+1)\theta(m),\label{thetam+kthetam}\\
\beta(m+k\theta(m))&=\frac{1}{k+1}\beta(m).\label{betam+kthetam}
\end{align}
\end{lem}
\begin{proof}
We prove this lemma by induction on $k\geq0$. The base case $k=0$ is trivial. Assume the case for $k\geq0$ has been proved,
that is, $m+k\theta(m)\in \mathcal{S}_\beta$ and Eq.~\eqref{thetam+kthetam} and \eqref{betam+kthetam} hold.
From the Rota-Baxter equation \eqref{rtequ} it follows that
$$P(x^m)P(x^{m+k\theta(m)})=P(x^mP(x^{m+k\theta(m)}))+P(P(x^m)x^{m+k\theta(m)}).$$ But now
\begin{align*}
P(x^m)P(x^{m+k\theta(m)})
=\beta(m)\beta(m+k\theta(m))x^{\theta(m)+\theta(m+k\theta(m))}=\frac{1}{k+1}\beta(m)^2x^{(k+2)\theta(m)},
\end{align*}
and
\begin{align*}
&P(x^mP(x^{m+k\theta(m)}))+P(P(x^m)x^{m+k\theta(m)})\\
=&\beta(m+k\theta(m))\beta(m+(k+1)\theta(m))x^{\theta(m+(k+1)\theta(m))}+\beta(m)\beta(m+(k+1)\theta(m))x^{\theta(m+(k+1)\theta(m))}\\
=&\frac{k+2}{k+1}\beta(m)\beta(m+(k+1)\theta(m))x^{\theta(m+(k+1)\theta(m))}.
\end{align*}
Since $\beta(m)\neq0$, we must have
\begin{align*}
\theta(m+(k+1)\theta(m))&=(k+2)\theta(m),\\
\beta(m+(k+1)\theta(m))&=\frac{1}{k+2}\beta(m).
\end{align*}
Clearly, $m+(k+1)\theta(m)\in \mathcal{S}_\beta$, which completes the induction.
\end{proof}

\begin{lem}\label{sbetasubset1}
Let $P$ be a monomial Rota-Baxter operator on $\mathbf{k}[x]$. Then
$\mathcal{S}_\beta+\theta(\mathcal{S}_\beta)\subseteq \mathcal{S}_\beta$.
\end{lem}
\begin{proof}
Assume that there exist $m,n\in \mathcal{S}_\beta$ such that $m+\theta(n)\not\in\mathcal{S}_\beta$.
By induction on $k\geq0$, we first show that $\theta(m)+n+k\theta(n)\in\mathcal{S}_\beta$, and
\begin{align}
\theta(\theta(m)+n+k\theta(n))&=\theta(m)+(k+1)\theta(n),\label{thetatheta(m)ktheta(n)n)}\\
\beta(\theta(m)+n+k\theta(n))&=\frac{1}{k+1}\beta(n).\label{betathetsfsdsdetfsn}
\end{align}
It follows from $m+\theta(n) \in\mathcal{Z}_\beta$ that $\beta(m+\theta(n))=0$.
Applying the definition of the Rota-Baxter operator $P$ gives that $P(x^m)P(x^n)=P(x^mP(x^n))+P(P(x^m)x^n)$, that is,
\begin{align*}
\beta(m)\beta(n)x^{\theta(m)+\theta(n)}&=\beta(m+\theta(n))\beta(n)x^{\theta(m+\theta(n))}+\beta(m)\beta(\theta(m)+n)x^{\theta(\theta(m)+n)}\\
&=\beta(m)\beta(\theta(m)+n)x^{\theta(\theta(m)+n)}.
\end{align*}
Notice that $m,n\in \mathcal{S}_\beta$ imply $\beta(m)\beta(n)\neq0$, so we have $\theta(\theta(m)+n)=\theta(m)+\theta(n)\neq0$  and $\beta(\theta(m)+n)=\beta(n)\neq0$,
whence $\theta(m)+n\in \mathcal{S}_\beta$. This shows Eq. \eqref{thetatheta(m)ktheta(n)n)} and \eqref{betathetsfsdsdetfsn}
hold for $k=0$.

Now assume that Eq. \eqref{thetatheta(m)ktheta(n)n)} and \eqref{betathetsfsdsdetfsn} is true for $k\geq0$. Considering the equation
\begin{align}\label{pmknkn1}
P(x^{\theta(m)+n+k\theta(n)})P(x^n)=P(x^{\theta(m)+n+k\theta(n)}P(x^n))+P(P(x^{\theta(m)+n+k\theta(n)})x^n).
\end{align}
According to Definition \ref{defnpbttht} and the inductive assumption,  the left hand side of Eq. \eqref{pmknkn1} is
\begin{align*}
\beta({\theta(m)+n+k\theta(n)})\beta(n)x^{\theta({\theta(m)+n+k\theta(n)})+\theta(n)}
=\frac{1}{k+1}\beta(n)^2x^{\theta(m)+(k+2)\theta(n)},
\end{align*}
the right hand side of Eq. \eqref{pmknkn1} is
\begin{align*}
&\left[\beta(n)\beta(\theta(m)+n+(k+1)\theta(n))+
\beta(\theta(m)+n+k\theta(n))\beta(\theta(\theta(m)+n+k\theta(n))+n)\right]x^{\theta(\theta(m)+n+(k+1)\theta(n))}\cr
=&\frac{k+2}{k+1}\beta(n)\beta(\theta(m)+n+(k+1)\theta(n))x^{\theta(\theta(m)+n+(k+1)\theta(n))}.
\end{align*}
Consequently, we have
\begin{align*}
\frac{1}{k+1}\beta(n)^2x^{\theta(m)+(k+2)\theta(n)}=\frac{k+2}{k+1}\beta(n)\beta(\theta(m)+n+(k+1)\theta(n))x^{\theta(\theta(m)+n+(k+1)\theta(n))}.
\end{align*}
Then combining this with the fact that $\beta(n)\neq0$ it follows that that
\begin{align*}
\theta(\theta(m)+n+(k+1)\theta(n))&=\theta(m)+(k+2)\theta(n),\\
\beta(\theta(m)+n+(k+1)\theta(n))&=\frac{1}{k+2}\beta(n),
\end{align*}
which completes the inductive argument. Therefore, Eq. \eqref{thetatheta(m)ktheta(n)n)} and \eqref{betathetsfsdsdetfsn}
hold for all $k\geq0$.

In view of Lemma \ref{m+thetam}, $n+k\theta(n)\in \mathcal{S}_\beta$ for any nonnegative integer $k$.
Now, by using identities \eqref{thetam+kthetam}$-$\eqref{betathetsfsdsdetfsn}, we have
\begin{align*}
P(x^m)P(x^{n+k\theta(n)})&=\beta(m)\beta(n+k\theta(n))x^{\theta(m)+\theta(n+k\theta(n))}\\
&=\frac{1}{k+1}\beta(m)\beta(n)x^{\theta(m)+(k+1)\theta(n)},   &(\hbox{by  Eq.~\eqref{thetam+kthetam} and \eqref{betam+kthetam}})
\end{align*}
\vspace{-3mm}
\begin{align*}
P(x^mP(x^{n+k\theta(n)}))&=\beta(m+\theta(n+k\theta(n)))\beta(n+k\theta(n))x^{\theta(m+\theta(n+k\theta(n)))}\\
&=\frac{1}{k+1}\beta(m+(k+1)\theta(n))\beta(n)x^{\theta(m+(k+1)\theta(n))},   &(\hbox{by Eq.~\eqref{thetam+kthetam} and \eqref{betam+kthetam}})
\end{align*}
and
\begin{align*}
P(P(x^m)x^{n+k\theta(n)})&=\beta(m)\beta(\theta(m)+n+k\theta(n))x^{\theta(\theta(m)+n+k\theta(n))}\\
&=\frac{1}{k+1}\beta(m)\beta(n)x^{\theta(m)+(k+1)\theta(n)}.   &(\hbox{by Eq.~\eqref{thetatheta(m)ktheta(n)n)} and \eqref{betathetsfsdsdetfsn}})
\end{align*}
Thus, $P(x^m)P(x^{n+k\theta(n)})=P(P(x^m)x^{n+k\theta(n)})$. Using Eq. \eqref{rtequ},  we have
\begin{align*}
P(x^m)P(x^{n+k\theta(n)})=P(x^mP(x^{n+k\theta(n)}))+P(P(x^m)x^{n+k\theta(n)}),
\end{align*}
so we conclude that $P(x^mP(x^{n+k\theta(n)}))=0$, that is, $\beta(m+(k+1)\theta(n))=0$ for any nonnegative integer $k$, whence
$m+(k+1)\theta(n)\in \mathcal{Z}_\beta$. In particular, we have $m+\theta(m)\theta(n)\in \mathcal{Z}_\beta$ since $\theta(m)$ is a positive integer.
On the other hand, by Lemma \ref{m+thetam}, $m+\theta(m)\theta(n)\in \mathcal{S}_\beta$ since  $\theta(n)$ is a positive integer.
This is a contradiction, proving $\mathcal{S}_\beta+\theta(\mathcal{S}_\beta)\subseteq \mathcal{S}_\beta$, as required.
\end{proof}

We also revised \cite[Proposition 3.18(2)]{rbopaiaoguo} as follows.
\begin{prop}\label{decompositionofsbeta}
If $P$ is a nonzero monomial Rota-Baxter operator on $\mathbf{k}[x]$, then there exists $k\in \mathbb{N}^*$ such that
$$
\mathcal{S}_\beta=(s_1+d\mathbb{N})\uplus(s_2+d\mathbb{N})\uplus\cdots\uplus(s_k+d\mathbb{N}),
$$
where $d$ is the greatest common divisor of $\theta({\mathcal{S}_\beta})$, and $0\leq s_1<s_2<\cdots<s_k\leq d-1$ are all integers.
\end{prop}
\begin{proof}
Denote $T=\theta(\mathcal{S}_\beta)$. It follows from Eq.~\eqref{thetaeq3} that $T$ is a subsemigroup of $\mathbb{N}$. Write $d=gcd(T)$.
Then $T/d:=\{\frac{t}{d}|t\in T\}$ is a numerical semigroups
\cite[Proposition 10.1]{RGFinitelyGeneratedCommutativeMonoids} or \cite[Lemma 2.1]{RosalesPAGarcia-SanchezNumericalsemigroups}, meaning a
subsemigroup of $\mathbb{N}$ with a finite complement $G\subseteq \mathbb{N}$  of so-called gaps.
Thus we obtain $T =dN\backslash G$. We write $f\in \mathbb{N}$ for the conductor of $T/d$,
which is the least integer $x$ such that $x+\mathbb{N} \subseteq T/d$. Then $fd+d\mathbb{N}\subseteq T$ holds.

Let $\Omega_i=\mathcal{S}_\beta\cap (i+d\mathbb{N})$ for $i\in\{0,1,\cdots,d-1\}$. Then $\mathcal{S}_\beta=\uplus_{i=0}^{d-1}\Omega_i$.
We claim that either $\Omega_i=i+d\mathbb{N}$ or $\Omega_i=\emptyset$. Indeed, suppose $m\in \mathcal{S}_\beta\cap (i+d\mathbb{N})$,
$n\in \mathcal{Z}_\beta\cap (i+d\mathbb{N})$. Then $d|(m-n)$. Since $P$ is a monomial Rota-Baxter operator on $\mathbf{k}[x]$,
according to Theorem \ref{equicodmonomial}\eqref{bicontain1}, we have
\begin{align*}
m+fd+d\mathbb{N}\subseteq \mathcal{S}_\beta+ \theta(\mathcal{S}_\beta)\subseteq \mathcal{S}_\beta,
\qquad
n+fd+d\mathbb{N}\subseteq  \mathcal{Z}_\beta+ \theta(\mathcal{S}_\beta)\subseteq \mathcal{Z}_\beta.
\end{align*}
Thus,  $(m+fd+d\mathbb{N})\cap (n+fd+d\mathbb{N})\subseteq \mathcal{S}_\beta \cap\mathcal{Z}_\beta=\emptyset$, contradicting $d|(m-n)$.
Suppressing the empty ones among $\Omega_i$ to conclude that
$$
\mathcal{S}_\beta=(s_1+d\mathbb{N})\uplus(s_2+d\mathbb{N})\uplus\cdots\uplus(s_k+d\mathbb{N})\quad{\rm with}\quad 0\leq s_1<s_2<\cdots<s_k\leq d-1,
$$
as required.
\end{proof}

An immediate consequence of Proposition \ref{decompositionofsbeta} is the following result about the values of $\theta$.
\begin{lem}\label{valuyetheta}
Let $n\in \mathbb{N}$ and let $n\equiv\overline{n}({\rm mod}\ d)$, where $d=gcd(\theta(\mathcal{S}_\beta))$
 and $\overline{n}\in\{0,1,\cdots,d-1\}$. Then $n\in \mathcal{S}_\beta$  if and only if
$\overline{n}\in \{s_1,s_2,\cdots,s_k\}$. Moreover,
$\theta(n)=\theta(\overline{n})+n-\overline{n}$ for any $n\in \mathcal{S}_\beta$.
\end{lem}
\begin{proof}
According to  Proposition  \ref{decompositionofsbeta}, $n\in \mathcal{S}_\beta$ is equivalent to $\overline{n}\in \{s_1,s_2,\cdots,s_k\}$.

For the second part, defining a map $\widetilde{\theta}: \mathcal{S}_\beta\rightarrow \mathbb{Z}$ by $\widetilde{\theta}(n)=\theta(n)-n$,
one obtains from Eq.~\eqref{thetaeq3} that
$\widetilde{\theta}(n)=\widetilde{\theta}(n+\theta(m))$ for all $m,n\in \mathcal{S}_\beta$. Thus $\widetilde\theta$ is periodic, and suppose
$e$ is the primitive period of $\widetilde\theta$. Clearly, every $\theta(m)$ is a period of $\widetilde\theta$, so $e|\theta(m)$
for all $m\in\mathcal{S}_\beta$, which implies that $e|d$. On the other hand,
$\theta(s_1+e)=\widetilde{\theta}(s_1+e)+s_1+e=\widetilde{\theta}(s_1)+s_1+e=\theta(s_1)+e$, so $e=\theta(s_1+e)-\theta(s_1)$ and hence $d|e$,
whence $e=d$ holds. Thus, $d$ is the primitive period of $\widetilde\theta$.
If we write $n=l_nd+\overline{n}$, then  $\theta(n)=\widetilde{\theta}(l_nd+\overline{n})+n=\widetilde{\theta}(\overline{n})+n=\theta(\overline{n})+n-\overline{n}$, as required.
\end{proof}

We next give a formula for the values of $\beta$.

\begin{lem}\label{valubgseta}
Let $n\in \mathcal{S}_\beta$ with $n\equiv\overline{n}({\rm mod}\ d)$, where $\overline{n}\in\{0,1,\cdots,d-1\}$.
Then
\begin{align*}%\label{thevad98esofbeta}
\beta(n)=\frac{\theta(\overline{n})}{\theta(\overline{n})+n-\overline{n}}\beta(\overline{n}).
\end{align*}
\end{lem}
\begin{proof}
Take $m,n\in \mathcal{S}_\beta$ such that $\overline{m}=\overline{n}$.
Then, by Lemma \ref{valuyetheta}, Eq.~\eqref{betaeq4} yields that
\begin{align*}
\beta(n)\beta(m)=(\beta(n)+\beta(m))\beta(\theta(\overline{n})+m+n-\overline{n}).
\end{align*}
In view of Proposition \ref{decompositionofsbeta}, $\theta(\overline{n})+m+n-\overline{n}\in \mathcal{S}_\beta$, and thus
\begin{align*}%\label{btamnm1}
\frac{1}{\beta(\theta(\overline{n})+m+n-\overline{n})}=\frac{1}{\beta(m)}+\frac{1}{\beta(n)}
\end{align*}
holds in the quotient field of $\mathbf{k}$. In particular, for any $m_1,m_2,n_1,n_2\in \mathcal{S}_\beta$ such that
$m_1+n_1=m_2+n_2$ and
$\overline{m_1}=\overline{m_2}=\overline{n_1}=\overline{n_2}$, we have
\begin{align*}%\label{btamnm2}
\frac{1}{\beta(m_1)}+\frac{1}{\beta(n_1)}=\frac{1}{\beta(m_2)}+\frac{1}{\beta(n_2)}.
\end{align*}
Therefore, for any element $n=l_nd+\overline{n}$ of $\mathcal{S}_\beta$ with $l_n\geq 1$, we have $\overline{n}\in \mathcal{S}_\beta$ and
\begin{align}\label{btamnm3}
\frac{1}{\beta(n)}&=\frac{1}{\beta((l_n-1)d+\overline{n})}+\frac{1}{\beta(d+\overline{n})}-\frac{1}{\beta(\overline{n})}\cr
&=\frac{1}{\beta((l_n-2)d+\overline{n})}+\frac{2}{\beta(d+\overline{n})}-\frac{2}{\beta(\overline{n})}\cr
&=\cdots\ \cdots\cr
&=\frac{l_n}{\beta(d+\overline{n})}-\frac{l_n-1}{\beta(\overline{n})}.
\end{align}

Note that $d|\theta(\overline{n})$, so we may suppose $\theta(\overline{n})=\sigma_{\overline{n}}d$ for some positive integer $\sigma_{\overline{n}}$.
Then, by Eq.~\eqref{btamnm3}, we have
\begin{align}\label{btamnm4}
\frac{1}{\beta(\overline{n}+\theta(\overline{n}))}=\frac{\sigma_{\overline{n}}}{\beta(d+\overline{n})}-\frac{\sigma_{\overline{n}}-1}{\beta(\overline{n})}.
\end{align}
On the other hand, it follows from Eq.~\eqref{betam+kthetam} that
$\beta(\overline{n}+\theta(\overline{n}))=\frac{1}{2}\beta(\overline{n})$, which together with Eq.~\eqref{btamnm4} yields that
\begin{align*}
\frac{1}{\beta(d+\overline{n})}=\frac{\sigma_{\overline{n}}+1}{\sigma_{\overline{n}}}\frac{1}{\beta(\overline{n})},
\end{align*}
and hence, by Eq.~\eqref{btamnm3} again, we obatin
\begin{align}\label{btamnm5}
\beta(n)=\frac{\sigma_{\overline{n}}}{\sigma_{\overline{n}}+l_n}\beta(\overline{n})
=\frac{\theta(\overline{n})}{\theta(\overline{n})+n-\overline{n}}\beta(\overline{n}),
\end{align}
as required. Notice that if $n\leq d-1$, that is, $n=\overline{n}$, then the Eq.~\eqref{btamnm5} holds trivially. This completes the proof.
\end{proof}

Now we give a detailed classification for monomial Rota-Baxter operators $P$ of weight zero on $\mathbf{k}[x]$.

\begin{thm}\label{valbtheuyeaq2a}
Let $P$ be a monomial linear operator of weight $0$ on $\mathbf{k}[x]$ defined by $P(x^n)=\beta(n)x^{\theta(n)}$, $n\in \mathbb{N}$.
Then $P$ is a Rota-Baxter operator if and only if there exist a positive integer $d$;  $d$ nonnegative integers
$c_0,c_1,\cdots,c_{d-1}$;
and $d$ elements $b_0,b_1,\cdots,b_{d-1}\in \mathbf{k}$  such that
\begin{enumerate}
\item\label{weight01main} $b_i=0$ if and only if $c_i=0$ where $i=0,1,\cdots, d-1$;

\item\label{weight02main} for all $n\in\mathbb{N}$, we have
 \begin{equation}\label{maintheq01}
\theta(n)=
\begin{cases}
0,& b_{\overline{n}}=0,\\
c_{\overline{n}}d+n-\overline{n},& b_{\overline{n}}\neq0,
\end{cases}
\end{equation}
and
\begin{equation}\label{maintheq02}
\beta(n)=
\begin{cases}
0,& b_{\overline{n}}=0,\\
\frac{b_{\overline{n}}c_{\overline{n}}d}{c_{\overline{n}}d+n-\overline{n}},&  b_{\overline{n}}\neq0,
\end{cases}
\end{equation}
where $\overline{n}\in\{0,1,\cdots,d-1\}$ is the remainder of $n$ module $d$.
\end{enumerate}
\end{thm}

\begin{proof}
For the case of $P=0$, it is enough to take $d=1$, $b_0=0$ and $c_0=0$.
In what follows, we assume that $P$ is a nonzero operator.

It follows from $P\neq0$ that $\mathcal{S}_\beta\neq \emptyset$,
we let $d=gcd(\theta(\mathcal{S}_\beta))$,
$b_i=\beta(i)$ and $c_i=\frac{1}{d}\theta(i)$, where $i=0,1,\cdots,d-1$.
Then, $d$ is a positive integer, $c_i\in\mathbb{N}$, $b_i\in \mathbf{k}$ for all $i=0,1,\cdots,d-1$.
From the fact that $\mathcal{S}_\beta=\mathcal{S}_\theta$ and
$\mathcal{Z}_\beta=\mathcal{Z}_\theta$ we see that \eqref{weight01main}  holds, while Lemmas \ref{valuyetheta} and \ref{valubgseta}
guarantee Eq.~\eqref{maintheq01} and \eqref{maintheq02} hold, respectively.

To prove the converse, we only need to show that the $\theta$ and $\beta$ defined in the theorem satisfy the
conditions \eqref{bicontain1} and \eqref{bicontain2} of Theorem \ref{equicodmonomial}.
If $b_i=0$ for all $i=0,1,\cdots,d-1$, then $P=0$ is trivial. If $b_i$ are not all zero, then $\mathcal{S}_\beta\neq\emptyset$.
By conditions \eqref{weight01main} and  \eqref{weight02main}, we have
\begin{align}\label{sbtb}
\mathcal{S}_\beta=\mathcal{S}_\theta=\biguplus_{b_i\neq0\atop 0\leq i\leq d-1 }(i+d\mathbb{N}),
\end{align}
and $\mathcal{Z}_\beta=\mathcal{Z}_\theta=\mathbb{N}\backslash \mathcal{S}_\beta$.
It's clearly that $gcd(\theta(\mathcal{S}_\beta))=d$ so that $\theta(\mathcal{S}_\beta)\subseteq d\mathbb{N}$,
and hence both $\mathcal{Z}_\beta+\theta(\mathcal{S}_\beta)\subseteq \mathcal{Z}_\beta$
and $\mathcal{S}_\beta+\theta(\mathcal{S}_\beta)\subseteq \mathcal{S}_\beta$ hold.
This concludes the  condition \eqref{bicontain1} of Theorem \ref{equicodmonomial} is satisfied.
Next we show Theorem \ref{equicodmonomial}\eqref{bicontain2}  also holds. To this end, taking any $m,n\in \mathcal{S}_\beta$.
By Eq.~\eqref{sbtb}, there exist $i,j\in\mathcal{S}_\beta\cap \{0,1,\cdots,d-1\}$ such that
$m=l_md+i$ and $n=l_nd+j$ for some $l_m,l_n\in \mathbb{N}$.
Then, by Eq.~\eqref{maintheq01}, we have
\begin{align*}
\theta(m)+\theta(n)&=c_id+m-i+c_jd+n-j\\
&=(c_i+l_m+c_j+l_n)d\\
&=\theta(l_md+i+c_jd+l_nd)\\
&=\theta(m+\theta(n)).
\end{align*}
Similarly, we also have $\theta(m)+\theta(n)=\theta(\theta(m)+n).$
By Eq.~\eqref{maintheq01} and \eqref{maintheq02},
\begin{align*}
\beta(m)\beta(n)&=\frac{b_ic_id}{(c_id+m-i)}\frac{b_jc_jd}{(c_jd+n-j)}
=\frac{b_ib_jc_ic_j}{(c_i+l_m)(c_j+l_n)}
\end{align*}
and
\begin{align*}
\beta(m)\beta(\theta(m)+n)+\beta(n)\beta(m+\theta(n))
&=\frac{b_ic_i}{c_i+l_m}\frac{b_jc_j}{c_i+c_j+l_m+l_n}+\frac{b_jc_j}{c_j+l_n}\frac{b_ic_i}{c_i+c_j+l_m+l_n}\\
&=\frac{b_ib_jc_ic_j}{(c_i+l_m)(c_j+l_n)}.
\end{align*}
Therefore, $\beta(m)\beta(n)=\beta(m)\beta(\theta(m)+n)+\beta(n)\beta(m+\theta(n))$ also holds. This completes the proof.
\end{proof}

Theorem \ref{valbtheuyeaq2a} gives a complete classification for all monomial Rota-Baxter operators of weight zero on ${\mathbf k}[x]$.
Now we give some examples.

\begin{exa}

(1) Take $d=1$ and $b_0=c_0=0$, then one obtains $P(x^n)=0$, so $P$ is the zero Rota-Baxter operator.

(2) Take $b_0,b_1,\cdots,b_{d-1}$ as nonzero elements of $\mathbf{k}$, one obtains $P(x^n)\neq0$ for all $n\in\mathbb{N}$. This will happen if and only if $P$ is nondegenerate.

(3) Take $d=1$, $c_0=c\in \mathbb{N}$ to be a positive integer and $b_0=\frac{b}{c}\in\mathbf{k}$ a nonzero element, one obtains $P(x^n)=\frac{b}{n+c}x^{n+c}=b\int_0^xt^{n+c-1}dt$. This exactly the case of $P$ is injective in view of
\cite[Theorem 3.13]{rbopaiaoguo}.
If we further take $b=c=1$, then $P(x^n)=\frac{1}{n+1}x^{n+1}$, and $P$ is the standard integration operator.
\end{exa}

\section{Monomial Rota-Baxter operators of weight nonzero}\label{sectionweightnonzero}
In this section, we investigate the classification of Rota-Baxter operators on $\mathbf{k}[x]$ of weight nonzero.
All monomial Rota-Baxter operators $P$ are assumed to be of weight $\lambda\neq0$
defined by $P(x^n)=\beta(n)x^{\theta(n)}$, $n\in \mathbb{N}$.

We first give an example to point out that the cases of weight zero and  nonzero are different greatly from each other.
\begin{exa}\label{examp21}
Let $\lambda\in \mathbf{k}\backslash\{0\}$.
For all $n\in\mathbb{N}$, define $\theta:\mathbb{N}\rightarrow\mathbb{N}$ by $\theta(n)=0$,
and $\beta:\mathbb{N}\rightarrow \mathbf{k}$ by $\beta(n)=-\lambda$.
One can easily to check that $P:\mathbf{k}[x]\rightarrow \mathbf{k}[x]$ defined by $P(x^n)=\beta(n)x^{\theta(n)}=-\lambda $ is a monomial Rota-Baxter operator on $\mathbf{k}[x]$ of weight $\lambda$. Clearly, $\mathcal{S}_\beta=\mathcal{Z}_\theta=\mathbb{N}$, $\mathcal{S}_\theta=\mathcal{Z}_\beta=\emptyset$.
This is impossible for  monomial Rota-Baxter operators of weight zero, because $\mathcal{S}_\beta=\mathcal{S}_\theta$ and
$\mathcal{Z}_\beta=\mathcal{Z}_\theta$ for the case of weight zero.
\end{exa}

For convenience, we first give some identities for later use. Since $P$ is a monomial Rota-Baxter operator of weight $\lambda$ on $\mathbf{k}[x]$,
the Rota-Bxater relation in Eq.~\eqref{rtequ} is equivalent to
\begin{align*}
P(x^m)P(x^n)=P(x^mP(x^n))+P(P(x^m)x^n)+\lambda P(x^{m+n}),
\end{align*}
that is,
\begin{align}\label{lambdaneq01}
&\beta(m)\beta(n)x^{\theta(m)+\theta(n)}\cr
=&\beta(m+\theta(n))\beta(n)x^{\theta(m+\theta(n))}+\beta(m)\beta(\theta(m)+n)x^{\theta(\theta(m)+n)}+\lambda\beta(m+n)x^{\theta(m+n)}
\end{align}
holds for all $m,n\in \mathbb{N}$.
If all the coefficients in Eq.~\eqref{lambdaneq01} are nonzero, then we must have either all the exponents of $x$ are equal or two of them are equal and
the other two are equal. We will use this fact frequently but no explanation in this section.

Let $m=n$ in Eq.~\eqref{lambdaneq01}, one obtains
\begin{align}\label{lambdaneq02}
\beta(n)^2x^{2\theta(n)}=2\beta(n+\theta(n))\beta(n)x^{\theta(n+\theta(n))}+\lambda\beta(2n)x^{\theta(2n)}.
\end{align}
Taking $m=0$ in Eq.~\eqref{lambdaneq01}, we have
\begin{align}\label{lambdaneq03}
\beta(0)\beta(n)x^{\theta(0)+\theta(n)}
=\beta(\theta(n))\beta(n)x^{\theta(\theta(n))}+\beta(0)\beta(\theta(0)+n)x^{\theta(\theta(0)+n)}+\lambda\beta(n)x^{\theta(n)}.
\end{align}
Taking $m=n=0$, then Eq.~\eqref{lambdaneq01} yields that
\begin{align}\label{lambdaneq04}
\beta(0)^2x^{2\theta(0)}
=2\beta(0)\beta(\theta(0))x^{\theta(\theta(0))}+\lambda\beta(0)x^{\theta(0)}.
\end{align}

Next we give some properties about the mappings $\beta$ and $\theta$, which is critical for the main result.
\begin{lem}\label{nonzerobetaneq00}
Let $P$ be a monomial Rota-Baxter operator of weight $\lambda$ on $\mathbf{k}[x]$ defined by $P(x^n)=\beta(n)x^{\theta(n)}$, $n\in \mathbb{N}$.
Then
\begin{enumerate}
\item\label{non0proty1} $\theta(0)=0$ and $\theta(\theta(n))=\theta(n)$ for all $n\in\mathbb{N}$;

\item\label{non0proty2} either $\beta(\theta(n))=-\lambda$ for all $n\in\mathbb{N}$ or
\begin{equation*}%\label{thtathtan}
\beta(\theta(n))=
\begin{cases}
0,& n\in\mathcal{Z}_\beta,\\
-\lambda,& n\in\mathcal{S}_\beta.
\end{cases}
\end{equation*}
In particular, $\theta(\mathcal{S}_\beta)\subseteq \mathcal{S}_\beta$, $\beta(0)$ is either $0$ or $-\lambda$;

\item\label{non0proty3} $\mathcal{Z}_\theta$ is a subsemigroup of $\mathbb{N}$. If $\mathcal{Z}_\beta\neq \emptyset$, then $\mathcal{Z}_\theta=\{0\}\cup\mathcal{Z}_\beta$;

\item\label{non0proty4} $im(\theta)$ is a subsemigroup of $\mathbb{N}$, and $im(\theta)\cap\mathcal{Z}_\theta=\{0\}$.
\end{enumerate}
\end{lem}

\begin{proof}
 We complete the proof of \eqref{non0proty1} and \eqref{non0proty2} by considering the following two cases.

{\bf Case 1.} $\beta(0)=0$. Then $\theta(0)=0$ since $\mathcal{Z}_\beta\subseteq \mathcal{Z}_\theta$, so Eq.~\eqref{lambdaneq03} is equivalent to
\begin{align}\label{lambdaneq006}
\beta(\theta(n))\beta(n)x^{\theta(\theta(n))}+\lambda\beta(n)x^{\theta(n)}=0.
\end{align}
If $n\in\mathcal{S}_\beta$, then $\beta(n)\neq0$, so Eq.~\eqref{lambdaneq006} yields that $\beta(\theta(n))=-\lambda$ and $\theta(\theta(n))=\theta(n)$.
If $n\in\mathcal{Z}_\beta$, then, by Definition \ref{defnpbttht}, $n\in\mathcal{Z}_\theta$,
whence $\beta(\theta(n))=\beta(0)=0$ and $\theta(\theta(n))=\theta(0)=0=\theta(n)$, as required.

{\bf Case 2.} $\beta(0)\neq0$.
It follows from Eq.~\eqref{lambdaneq04} that
\begin{align}\label{lambdaneq05}
\beta(0)x^{2\theta(0)}=2\beta(\theta(0))x^{\theta(\theta(0))}+\lambda x^{\theta(0)}.
\end{align}
Consequently, $2\theta(0)=\theta(0)$ so that $\theta(0)=0$. Thus, by Eq.~\eqref{lambdaneq05}, $\beta(0)=-\lambda$.
By Eq.~\eqref{lambdaneq03} again, we can also get Eq.~\eqref{lambdaneq006}.
Then, for any $n\in\mathcal{S}_\beta$, by Eq.~\eqref{lambdaneq006},
$\beta(\theta(n))=-\lambda$ and $\theta(\theta(n))=\theta(n)$.
If $n\in\mathcal{Z}_\beta$, then $n\in\mathcal{Z}_\theta$,
whence $\beta(\theta(n))=\beta(0)=-\lambda$ and $\theta(\theta(n))=\theta(0)=0=\theta(n)$, as required.

\eqref{non0proty3} We prove the desired results via proving $\mathcal{Z}_\beta$, $\mathcal{Z}_\theta$ and $\mathcal{Z}_\theta\cap \mathcal{S}_\beta$ are
all subsemigroups of $\mathbb{N}$ if they are nonempty.

Let $m\in\mathcal{Z}_\beta$ and $n\in\mathcal{Z}_\theta$. It follows from Eq.~\eqref{lambdaneq01} that $\lambda\beta(m+n)x^{\theta(m+n)}=0$.
So $\lambda\neq0$ yields $\beta(m+n)=0$, and hence $m+n\in\mathcal{Z}_\beta$ so that $\mathcal{Z}_\beta+\mathcal{Z}_\theta\subseteq\mathcal{Z}_\beta$.
In particular, $\mathcal{Z}_\beta$ is a subsemigroup of $\mathbb{N}$ since $\mathcal{Z}_\beta\subseteq\mathcal{Z}_\theta$.

If we suppose that $m,n\in\mathcal{Z}_\theta$, then Eq.~\eqref{lambdaneq01} is equivalent to
\begin{align}\label{betammnmqw10}
\beta(m)\beta(n)+\lambda\beta(m+n)x^{\theta(m+n)}=0.
\end{align}
If at least one of $m,n$, say $m$, in $\mathcal{Z}_\beta$, then
$m+n\in\mathcal{Z}_\beta+\mathcal{Z}_\theta\subseteq\mathcal{Z}_\beta$ as has been proved so that
$m+n\in\mathcal{Z}_\theta$.
If $m,n\in \mathcal{Z}_\theta\backslash\mathcal{Z}_\beta$, then $\beta(m)\beta(n)\neq0$ so that $\theta(m+n)=0$ by Eq.~\eqref{betammnmqw10}. Thus,
 we also have $m+n\in \mathcal{Z}_\theta$.
This shows $\mathcal{Z}_\theta$ is a subsemigroup of $\mathbb{N}$.

Taking any $m,n\in \mathcal{Z}_\theta\cap\mathcal{S}_\beta$, then Eq.~\eqref{lambdaneq01} is equivalent to
$\beta(m)\beta(n)+\lambda\beta(m+n)x^{\theta(m+n)}=0$.
Notice that $\lambda$, $\beta(m)$ and $\beta(n)$ are all nonzero,
so $\beta(m+n)\neq0$ and $\theta(m+n)=0$, that is, $m+n\in\mathcal{S}_\beta\cap \mathcal{Z}_\theta$.
Therefore, $\mathcal{Z}_\theta\cap \mathcal{S}_\beta$ is also a subsemigroup of $\mathbb{N}$.

Now, assume that $\mathcal{Z}_\beta\neq\emptyset$. By \eqref{non0proty1}, $0\in \mathcal{Z}_\theta$ holds.
If $\mathcal{Z}_\beta=\mathcal{Z}_\theta$, then $0\in \mathcal{Z}_\beta$, and hence $\mathcal{Z}_\theta=\{0\}\cup\mathcal{Z}_\beta$.
If $\mathcal{Z}_\beta\neq\mathcal{Z}_\theta$, then $\mathcal{Z}_\theta\cap\mathcal{S}_\beta\neq \emptyset$, and
the subsemigroup $\mathcal{Z}_\theta$ is a disjoint union of $\mathcal{Z}_\beta$ and $\mathcal{Z}_\theta\cap\mathcal{S}_\beta$,
which are also two subsemigroups of $\mathbb{N}$, so one of $\mathcal{Z}_\beta$ and $\mathcal{Z}_\theta\cap\mathcal{S}_\beta$ must be $\{0\}$.
Notice that $\mathcal{Z}_\beta+\mathcal{Z}_\theta\subseteq\mathcal{Z}_\beta$, so $\mathcal{Z}_\theta\cap\mathcal{S}_\beta=0$.
 Thus, in either case, we must have
$\mathcal{Z}_\theta=\{0\}\cup\mathcal{Z}_\beta$ holds.

\eqref{non0proty4} Taking $s,t\in im(\theta)$. Without loss of generality, assume that $s,t\neq0$.
Then there exist $m,n\in \mathcal{S}_\theta$ such that $s=\theta(m)$ and $t=\theta(n)$.
Since $\mathcal{S}_\theta\subseteq \mathcal{S}_\beta$, we have $m,n\in \mathcal{S}_\beta$.
 By Eq.~\eqref{lambdaneq01},
$$
s+t=\theta(m)+\theta(n)\in\{\theta(m+\theta(n)),\theta(\theta(m)+n),\theta(m+n)\}\subseteq im(\theta).
$$
So $im(\theta)$ is a subsemigroup of $\mathbb{N}$.

Let $m\in im(\theta)\cap\mathcal{Z}_\theta$. Then, by \eqref{non0proty1}, $m=\theta(m)=0$ so that $im(\theta)\cap\mathcal{Z}_\theta\subseteq\{0\}$.
On the other hand, $\theta(0)=0$ yields that
$0\in im(\theta)\cap\mathcal{Z}_\theta$ and hence $im(\theta)\cap\mathcal{Z}_\theta=\{0\}$.
\end{proof}

\begin{lem}\label{m=2k+1n+theta(n)}
Let $P$ be a nonzero monomial Rota-Baxter operator of weight $\lambda$ on $\mathbf{k}[x]$ defined by $P(x^n)=\beta(n)x^{\theta(n)}$, $n\in \mathbb{N}$,
where $\mathcal{S}_\theta=\mathbb{N}^*$. Then
\begin{enumerate}
\item\label{nn2n2kn} for any $k,n\in \mathbb{N}^*$,  $\theta(n+\theta(n))=2\theta(n)$ and $\theta(2^kn)=2^k\theta(n)$;

\item\label{mnmnm2mn} for any $m,n\in \mathbb{N}^*$, $\theta(m)=\theta(n)$ implies $\theta(m+n)=2\theta(m)$.
\end{enumerate}
\end{lem}
\begin{proof}
\eqref{nn2n2kn} By Eq.~\eqref{lambdaneq02}, we have
\begin{align*}%\label{m=neq}
\theta(2n)=\theta(n+\theta(n))=2\theta(n).
\end{align*}
Clearly, one has $\theta(2^kn)=2\theta(2^{k-1}n)=\cdots=2^k\theta(n)$.

\eqref{mnmnm2mn} Suppose that $\theta(m)=\theta(n)$. Then one has
$\theta(m+\theta(n))=\theta(m+\theta(m))=2\theta(m)$ by \eqref{nn2n2kn}. By symmetry, $\theta(\theta(m)+n)=2\theta(n)$ and hence $\theta(\theta(m)+n)=\theta(m+\theta(n))$.
Then, it follows from Eq.~\eqref{lambdaneq01} that $\theta(m+n)=\theta(m)+\theta(n)=2\theta(m)$, as required.
\end{proof}

\begin{lem}\label{Sbeta=nNn}
Let $P$ be a nonzero monomial Rota-Baxter operator of weight $\lambda$ on $\mathbf{k}[x]$ defined by $P(x^n)=\beta(n)x^{\theta(n)}$, $n\in \mathbb{N}$,
where $\mathcal{S}_\theta=\mathbb{N}^*$.
Then $\theta(n)=n$ for all $n\in\mathbb{N}$.
\end{lem}
\begin{proof}
Since $\mathcal{S}_\theta\subseteq\mathcal{S}_\beta$, one has $\mathbb{N}^*\subseteq\mathcal{S}_\beta$.
It follows from $\mathcal{S}_\theta=\mathbb{N}^*$ that $im(\theta)\neq\{0\}$.
Let $d=gcd(im(\theta))$, and then $\frac{1}{d}im(\theta)$ is a numerical semigroup by Lemma \ref{nonzerobetaneq00}\eqref{non0proty4}.
We write $f$ for the conductor of $\frac{1}{d}im(\theta)$.
Then $df+d\mathbb{N}\subseteq im(\theta)$. In particular, for an enough large $k\in \mathbb{N}$, we must have $2^kd\in im(\theta)$.
Thus, by Lemma \ref{nonzerobetaneq00}\eqref{non0proty1} and Lemma \ref{m=2k+1n+theta(n)}\eqref{nn2n2kn}, $2^kd=\theta(2^kd)=2^k\theta(d)$ so that $\theta(d)=d$, which means that
$d=gcd(im(\theta))\in im(\theta)$ and hence $im(\theta)=d\mathbb{N}$.

Clearly, $d$ must be an odd number. Otherwise, $\frac{d}{2}=\frac{1}{2}\theta(d)=\theta(\frac{d}{2})\in im(\theta)=d\mathbb{N}$, a contradiction.
We claim that $d=1$. Assume the contrary $d\geq 3$ holds. Take any two positive integers $m,n$ such that $m+n=d$, then one has
$\theta(m+n)=d$. Notice that $\theta(m),\theta(n)\in d\mathbb{N}^*$, so $\theta(m)+\theta(n)\neq d=\theta(m+n)$.
By Eq.~\eqref{lambdaneq01}, without loss of generality, suppose that
\begin{align}
\theta(m+\theta(n))&=\theta(m+n)=d,  \label{eq2.11}\\
\theta(\theta(m)+n)&=\theta(m)+\theta(n).\label{eq2.12}
\end{align}

We now proceed to obtain a contradiction via the following four steps:

{\bf Step 1.} We prove the following two identities by induction on $k$,
\begin{align}
\theta(n+kd)&=\theta(n)+kd,\label{eq2.21}\\
\theta(m+\theta(n)+kd)&=\theta(m+\theta(n))+kd \label{eq2.22}
\end{align}
for all $ k\in\mathbb{N}$.

Eq.~\eqref{eq2.21} and \eqref{eq2.22} are trivial for the case of $k=0$. Assume that Eq.~\eqref{eq2.21} and \eqref{eq2.22}
have been proved for $k\geq 0$. Replace $m,n$ by $m+\theta(n)$ and $n+kd$ in Eq.~\eqref{lambdaneq01} respectively, one has
$$
\begin{array}{llll}
\theta(m+\theta(n))+\theta(n+kd)&=d+\theta(n)+kd &(\hbox{by  Eq.~\eqref{eq2.11} and the induction hypothesis Eq.~\eqref{eq2.21}}) \cr
&=\theta(n)+(k+1)d,\cr
\theta(m+\theta(n)+\theta(n+kd))&=\theta(m+2\theta(n)+kd),&(\hbox{by the induction hypothesis Eq.~\eqref{eq2.21}}) \cr
\theta(\theta(m+\theta(n))+n+kd)&=\theta(n+(k+1)d),&(\hbox{by  Eq.~\eqref{eq2.11}}) \cr
\theta(m+\theta(n)+n+kd)&=\theta(\theta(n)+(k+1)d)\cr
&=\theta(n)+(k+1)d.&(\hbox{by Lemma \ref{nonzerobetaneq00}\eqref{non0proty1} and $im(\theta)=d\mathbb{N}$})
\end{array}
$$
Notice that $\mathbb{N}^*\subseteq \mathcal{S}_\beta$, so all the coefficients in Eq.~\eqref{lambdaneq01} are nonzero.
Comparing the exponents of $x$ in Eq.~\eqref{lambdaneq01} which are listed as above, we have
$$
\theta(m+2\theta(n)+kd)=\theta(n+(k+1)d).
$$
Then, by Lemma \ref{m=2k+1n+theta(n)}, one has
\begin{align*}
\begin{array}{llll}
\theta(n+(k+1)d)&=\frac{1}{2}\theta(m+2\theta(n)+kd+n+(k+1)d)&(\hbox{by  Lemma \ref{m=2k+1n+theta(n)}\eqref{mnmnm2mn}}) \cr
&=\frac{1}{2}\theta(2\theta(n)+2(k+1)d)&\cr
&=\theta(\theta(n)+(k+1)d) &(\hbox{by  Lemma \ref{m=2k+1n+theta(n)}\eqref{nn2n2kn}})\cr
&=\theta(n)+(k+1)d. &(\hbox{by $im\theta=d\mathbb{N}$  and Lemma \ref{nonzerobetaneq00}\eqref{nn2n2kn}})
\end{array}
\end{align*}
The induction hypothesis then yields Eq.~\eqref{eq2.21} holds for all $k\in \mathbb{N}$.

For Eq.~\eqref{eq2.22}, we substitute $n+kd$ for $n$ in Eq.~\eqref{lambdaneq01}. It follows from  Eq.~\eqref{eq2.21} that
$$
\begin{array}{llll}
\theta(m)+\theta(n+kd)&=\theta(m)+\theta(n)+kd, \cr
\theta(m+\theta(n+kd))&=\theta(m+\theta(n)+kd),\cr
\theta(\theta(m)+n+kd)&=\theta(m)+\theta(n)+kd,&(\hbox{by  Eq.~\eqref{eq2.21} and $\theta(m)+kd\in im(\theta)=d\mathbb{N}$})\cr
\theta(m+n+kd))&=\theta((k+1)d)=(k+1)d. &(\hbox{by $m+n=d$})
\end{array}
$$
Since all the coefficients in Eq\eqref{lambdaneq01} are nonzero, we have
$$\theta(m+\theta(n)+kd)=(k+1)d=\theta(m+\theta(n))+kd.$$ This yields Eq.~\eqref{eq2.22} holds.

{\bf Step 2.} Let $u=\frac{1}{d}max\{\theta(0),\theta(1),\cdots,\theta(d-1)\}$. We show that
\begin{align}\label{eq2.31}
\theta(s+ud+kd)=\theta(s+ud)+kd
\end{align}
for all $ k\in\mathbb{N}$ and all $s\in\{0,1,\cdots,d-1\}$.

Let $\theta(n)=ld$. Clearly, $l\leq u$ since $n\in\{0,1,\cdots,d-1\}$.
Then, by Eq.~\eqref{eq2.21} and \eqref{eq2.22},
\begin{align*}
\theta(n+ud+kd)&=\theta(n)+ud+kd=\theta(n+ud)+kd,\cr
\theta(m+ud+kd)&=\theta(m+\theta(n)+(u-l+k)d)\cr
&=\theta(m+\theta(n))+(u-l)d+kd\cr
&=\theta(m+\theta(n)+(u-l)d)+kd\cr
&=\theta(m+ud)+kd.
\end{align*}
By the arbitrariness of $m,n\in\mathbb{N}$ with $m+n=d$,
we obtain the desired result.

{\bf Step 3.} Let $\theta(1+ud)=cd$, we prove
\begin{align}\label{eq2.41}
\theta(s+ud)=scd-(s-1)ud \quad for\ all\quad s\in\{0,1,\cdots,d-1\}.
\end{align}

For $s=0$, $\theta(ud)=ud$ is clearly hold because of Lemma \ref{nonzerobetaneq00}\eqref{non0proty1} and $im(\theta)=d\mathbb{N}$. Assume that Eq.~\eqref{eq2.41} has been proved for $0\leq s\leq d-2$.
Then, by Eq.~\eqref{eq2.31}, one has
$$
\theta(1+ud+\theta(s+ud))=\theta(1+ud)+\theta(s+ud)=\theta(\theta(1+ud)+s+ud).
$$
Take $m=1+ud$ and $n=s+ud$ in Eq.~\eqref{lambdaneq01}, it follows that $\theta(1+ud+s+ud)=\theta(1+ud)+\theta(s+ud)$.
Note that $s+1\leq d-1$, then, by Eq.~\eqref{eq2.31} and the induction hypothesis,
we have
$$
\theta(s+1+ud)=\theta(1+ud+s+ud)-ud=\theta(1+ud)+\theta(s+ud)-ud=(s+1)cd-sud.
$$
So Eq.~\eqref{eq2.41} holds.

{\bf Step 4.} Getting a contradiction. Replace $m,n$ by $m'=1+ud$ and $n'=d-1 +ud$ in Eq.~\eqref{lambdaneq01}, respectively.
Then, by Eq.~\eqref{eq2.31},
one obtains that
$$
\theta(m')+\theta(n')=\theta(m'+\theta(n'))=\theta(\theta(m')+n'),
$$
whence
$\theta(m'+n')=\theta(m')+\theta(n')$. Since $m'+n'\in d\mathbb{N}=im(\theta)$, we have
$\theta(m'+n')=m'+n'=d+2ud$. By Lemma \ref{nonzerobetaneq00}\eqref{non0proty1} and Eq.~\eqref{eq2.41}, one has
$\theta(m')+\theta(n')=cd^2-(d-2)ud$,
that is,
$$
d+2ud=cd^2-(d-2)ud,
$$
and hence $(c-u)d=1$, contradicting $d\geq 3$.

Therefore, we must have $d=1$, and hence $im(\theta)=\mathbb{N}$.
In view of Lemma \ref{nonzerobetaneq00}\eqref{non0proty1}, $\theta(n)=n$ for all $n\in \mathbb{N}$.
\end{proof}

Now we establish the classification for monomial Rota-Baxter operators on $\mathbf{k}[x]$ of weight nonzero.

\begin{thm}\label{thmnonzero}
Let $P$ be a nonzero  monomial linear operator on $\mathbf{k}[x]$ of weight $\lambda\neq0$.
Then $P$ is a Rota-Baxter operator if and only if $P$ is one of the following cases:

\begin{enumerate}
\item\label{weightnzero1}  there exists $b\in \mathbf{k}\backslash\{0\}$ such that $P(x^n)=(-\lambda)^{1-n}b^n$ for all $n\in\mathbb{N}$;

\item\label{weightnzero2}   $P(x^n)=-\lambda x^n$ for all $n\in\mathbb{N}$;

\item\label{weightnzero3}   for all $n\in\mathbb{N}$,
\begin{equation*}
P(x^n)=
\begin{cases}
0,& n=0,\\
-\lambda x^n,& n\neq0;
\end{cases}
\end{equation*}

\item\label{weightnzero4}  for all $n\in\mathbb{N}$,
\begin{equation*}
P(x^n)=
\begin{cases}
-\lambda,& n=0,\\
0,& n\neq0.
\end{cases}
\end{equation*}
\end{enumerate}
\end{thm}

\begin{proof}
It is a routine to check that all the operators defined in \eqref{weightnzero1}-\eqref{weightnzero4} are monomial Rota-Baxter operators on $\mathbf{k}[x]$.
Conversely, let $P$ be a nonzero monomial Rota-Baxter operator of weight $\lambda$ on $\mathbf{k}[x]$
defined by $P(x^n)=\beta(n)x^{\theta(n)}$, $n\in \mathbb{N}$.
Now we prove $P$ must be one of the four types via the following cases.

{\bf Case 1.} $\mathcal{Z}_\beta=\emptyset$.

{\bf Case 1.1.}  $\mathcal{Z}_\theta\neq\{0\}$. By Lemma \ref{nonzerobetaneq00}\eqref{non0proty4},
the intersection of the two subsemigroups $im(\theta)$ and $\mathcal{Z}_\theta$ is $\{0\}$,
which means that $im(\theta)=\{0\}$, and thus
$\mathcal{Z}_\theta=\mathbb{N}$. Note that  $\mathcal{S}_\beta=\mathbb{N}$,
so Eq.~\eqref{lambdaneq01} is equivalent to
$\beta(m)\beta(n)+\lambda\beta(m+n)=0$.
Thus, for any $m_1,m_2,n_1,n_2\in \mathbb{N}$ such that $m_1+n_1=m_2+n_2$, we must have
$\beta(m_1)\beta(n_1)=\beta(m_2)\beta(n_2)$.
Notice that, in view of Lemma \ref{nonzerobetaneq00}\eqref{non0proty2}, $\beta(0)=-\lambda$, and then
it is easy to see that
\begin{align*}%\label{lambdaneq01e1}
\beta(n)=\frac{\beta(1)}{\beta(0)}\beta(n-1)=\left(\frac{\beta(1)}{\beta(0)}\right)^2\beta(n-2)=\cdots=\frac{\beta(1)^n}{\beta(0)^{n-1}}=(-\lambda)^{1-n}\beta(1)^n.
\end{align*}
Let $\beta(1)=b$, then $b\neq0$ and one has  $\beta(n)=(-\lambda)^{1-n}b^n$ for all $n\in\mathbb{N}$.
This case is reduced to \eqref{weightnzero1}.

{\bf Case 1.2.} $\mathcal{Z}_\theta=\{0\}$. Then $\mathcal{S}_\theta=\mathbb{N}^*$ and hence
$\theta(n)=n$ for all $n\in\mathbb{N}$ by Lemma \ref{Sbeta=nNn}.
It follows from Lemma \ref{nonzerobetaneq00}\eqref{non0proty2} that $\beta(n)=-\lambda$ for all $n\in\mathbb{N}$.
This case is reduced to \eqref{weightnzero2}.

{\bf Case 2.} $\mathcal{Z}_\beta\neq\emptyset$. By Lemma \ref{nonzerobetaneq00}\eqref{non0proty3}, $\mathcal{Z}_\theta=\{0\}\cup\mathcal{Z}_\beta$.

{\bf Case 2.1.} $\mathcal{Z}_\beta=\mathcal{Z}_\theta$. Then $0\in\mathcal{Z}_\theta=\mathcal{Z}_\beta$.

{\bf Case  2.1.1.} $\mathcal{Z}_\beta=\{0\}$.
Then $\mathcal{Z}_\theta=\{0\}$ and hence $\mathcal{S}_\beta=\mathcal{S}_\theta=\mathbb{N}^*$.
 Then, by Lemma \ref{Sbeta=nNn}, $\theta(n)=n$ for all $n\in\mathbb{N}$.
It follows from Lemma \ref{nonzerobetaneq00}\eqref{non0proty2} that $\beta(0)=0$ and $\beta(n)=-\lambda$ for all $n\in\mathbb{N}^*$.
This case is reduced to \eqref{weightnzero3}.

{\bf Case  2.1.2.} $\mathcal{Z}_\beta\neq\{0\}$. In this case $\mathcal{Z}_\beta=\mathcal{Z}_\theta$
are nonzero subsemigroup of $\mathbb{N}$. It follows from Lemma \ref{nonzerobetaneq00}\eqref{non0proty4} that $im(\theta)\cap\mathcal{Z}_\theta=\{0\}$,
whence $im(\theta)=\{0\}$, so that $\mathcal{Z}_\beta=\mathcal{Z}_\theta=\mathbb{N}$. This forces $P=0$, which contradicts $P$ is nonzero.

{\bf Case 2.2.} $\mathcal{Z}_\beta\neq\mathcal{Z}_\theta$.  By Lemma \ref{nonzerobetaneq00}\eqref{non0proty3}, one obtains that
$0\not\in\mathcal{Z}_\beta$ and $\mathcal{Z}_\theta=\{0\}\cup\mathcal{Z}_\beta$ so that $\mathcal{Z}_\beta$ is a nonzero subsemigroup of $\mathbb{N}$.
According to Lemma \ref{nonzerobetaneq00}\eqref{non0proty4},
$im\theta\cap\mathcal{Z}_\beta=\emptyset$. But $im(\theta)$ and $\mathcal{Z}_\beta$ both are subsemigroups of $\mathbb{N}$,
so $im(\theta)=\{0\}$. Therefore, we have
$\mathcal{Z}_\theta=\mathbb{N}$ whence $\mathcal{Z}_\beta=\mathbb{N}^*$, so
$\theta(n)=0$ for all $n\in\mathbb{N}$; by Lemma \ref{nonzerobetaneq00}\eqref{non0proty2}, $\beta(0)=-\lambda$ and $\beta(n)=0$ for all $n\in \mathbb{N}^*$.
This case is reduced to \eqref{weightnzero4}.
\end{proof}

The Rota-Baxter operators given by Theorem \ref{thmnonzero}\eqref{weightnzero1} and \eqref{weightnzero2} are nondegenerate, while those given by
Theorem \ref{thmnonzero}\eqref{weightnzero3} and \eqref{weightnzero4} are degenerate.  For a given $\lambda\neq0$ in $\mathbf{k}$, the addition of the operators defined by Theorem \ref{thmnonzero}\eqref{weightnzero3} and \eqref{weightnzero4} respectively gives the one defined by
Theorem \ref{thmnonzero}\eqref{weightnzero2}.
\begin{exa}
(1) For a given $\lambda\neq0$ in $\mathbf{k}$, put $b=-\lambda$. According to  Theorem \ref{thmnonzero}\eqref{weightnzero1}, the $\mathbf{k}$-linear
map $P:\mathbf{k}[x]\rightarrow\mathbf{k}[x]$ defined by $P(x^n)=-\lambda$ is a Rota-Baxter operator of weight $\lambda$. In this case we have $P(f(x))=-\lambda f(1)$
for any $f(x)\in \mathbf{k}[x]$.

(2) If we take $\lambda=-1$ and $b\in\mathbf{k}\backslash\{0\}$, then, by Theorem \ref{thmnonzero}\eqref{weightnzero1}, the operator $P$ defined by
$P(x^n)=b^n$ is Rota-Baxter. Moreover, we have $P(f(x))=f(b)$ for any $f(x)\in \mathbf{k}[x]$.

(3) The Rota-Baxter operator given by Theorem \ref{thmnonzero}\eqref{weightnzero2} is a scalar product. In particular, the identity map is a
Rota-Baxter operator on $\mathbf{k}[x]$ of weight $-1$.
\end{exa}

\smallskip

\noindent {\bf Acknowledgements}: The author would like to thank Professor Li Guo for the
valuable discussion and to Professor Gaoyong Zhang for his constant encouragement and support.
He would also like to thank the  NYU Polytechnic School of Engineering
for the excellent working conditions and to the China Scholarship Council
for the financial support that enabled he to go abroad.
This work was supported by NSFC grant $11426183$, NSF grant DMS $1312181$,
Chongqing Research Program of Application Foundation and Advanced Technology $($No. cstc2014jcyjA00028$)$ and
Fundamental Research Funds for the Central Universities $($No. XDJK2013C060$)$.

\end{document}